\documentclass{elsart3-1}

\usepackage{amssymb}

\usepackage[english,francais]{babel}

\newtheorem{theorem}{Theorem}[section]
\newtheorem{lemma}[theorem]{Lemma}
\newtheorem{e-proposition}[theorem]{Proposition}
\newtheorem{corollary}[theorem]{Corollary}
\newtheorem{e-definition}[theorem]{Definition\rm}
\newtheorem{remark}{\it Remark\/}

\newtheorem{theoreme}{Th\'eor\`eme}[section]

\newtheorem{proposition}[theoreme]{Proposition}

\renewcommand{\theequation}{\arabic{equation}}
\def\og{\leavevmode\raise.3ex\hbox{$\scriptscriptstyle\langle\!\langle$~}}
\def\fg{\leavevmode\raise.3ex\hbox{~$\!\scriptscriptstyle\,\rangle\!\rangle$}}

\journal{the Acad\'emie des sciences}
\begin{document}
\centerline{}
\begin{frontmatter}


\selectlanguage{english}
\title{Cram\'{e}r type  moderate deviations for stationary sequences of bounded random variables}


\selectlanguage{english}
\author{Xiequan Fan}
\ead{fanxiequan@hotmail.com}
\address{Center for Applied Mathematics,
Tianjin University, Tianjin, China}


\medskip
\begin{center}
{\small Received *****; accepted after revision +++++\\
Presented by £££££}
\end{center}

\begin{abstract}
\selectlanguage{english}
We derive  Cram\'{e}r type moderate deviations for stationary sequences of bounded random variables.
Our results imply the moderate deviation principles and a Berry-Esseen bound.  Applications to quantile coupling inequalities, functions of $\phi$-mixing sequences, and contracting Markov chains  are discussed.
{\it To cite this article: A.
Name1, A. Name2, C. R. Acad. Sci. Paris, Ser. I 340 (2005).}

\vskip 0.5\baselineskip

\selectlanguage{francais}
\noindent{\bf R\'esum\'e} \vskip 0.5\baselineskip \noindent
{\bf D\'{e}viations mod\'{e}r\'{e}s de type Cram\'{e}r  pour les s\'{e}quences stationnaires.}
Nous d\'{e}rivons les d\'{e}viations mod\'{e}r\'{e}es de type Cram\'{e}r pour des s\'{e}quences stationnaires de variables al\'{e}atoires born\'{e}es. Nos r\'{e}sultats impliquent les principes de d\'{e}viation mod\'{e}r\'{e}e et un th\'{e}oreme  de Berry-Esseen. Les applications aux  in\'{e}galit\'{e}s de couplage quantile, fonctions des s\'{e}quences de m\'{e}lange, et des cha\^{\i}nes de Markov contractantes   sont discut\'{e}es.
{\it Pour citer cet article~: A. Name1, A. Name2, C. R. Acad. Sci.
Paris, Ser. I 340 (2005).}
\end{abstract}
\end{frontmatter}

\section{Introduction}
For the stationary sequence $(X_i)_{i \in \mathbf{Z}}$ of centered random variables, define the partial sums and the normalized partial sums process by
\[
S_n=\sum_{i=1}^n X_i  \ \ \ \ \ \ \textrm{and} \ \ \ \ \   \ W_n = \frac1{\sqrt{n}}S_n,
\]
respectively. 
We say that the sequence of random variables $\{W_n, n >1\}$ satisfies the moderate deviation principle (MDP)   with speed $a_n \rightarrow 0$
and good rate function $I(\cdot)$, if the level set $\{x, \, I(x) \leq t \}$ are compact for all $t \in \mathbb{R},$ and for all
Borel sets $B,$
\begin{eqnarray}
- \inf_{x \in B^o}  I(x) &\leq & \liminf_{n\rightarrow \infty}   a_n^2  \ln \mathbf{P}\big(a_n  W_n     \in B \big)
  \leq  \limsup_{n\rightarrow \infty} a_n^2 \ln \mathbf{P}\big(a_n  W_n   \in B \big) \leq  - \inf_{x \in \overline{B}}I(x) ,
\end{eqnarray}
where $B^o$ denotes the interior of $B,$   $\overline{B}$  the closure of $B,$ and the infimum of a function over an empty set is interpreted as $\infty.$
 The MDP is an intermediate behavior between the central limit theorem ($a_n=O(1))$ and large deviations $(a_n \asymp \frac1{\sqrt{n}} ).$

The MDP results have been obtained by several authors. De Acosta  \cite{DA92} applied Laplace approximations to prove the MDP for sums of independent random vectors.   Dembo \cite{D96} showed that the MDP holds for the trajectory of a locally square integrable martingale with bounded jumps as soon as its quadratic covariation converges in  probability at an exponential rate.   Gao \cite{G96} and Djellout \cite{D02} obtained the MDP for martingales with non-bounded differences and $\phi$-mixing sequences with summable mixing rate. Dedecker et al.\ \cite{DMPU09} derived the MDP for stationary sequences of bounded random variables under martingale-type conditions. It is known that the MDP results for stationary sequences can be applied in a variety of settings. For instance,  Dedecker et al.\ \cite{DMPU09} showed that such type of results can be applied to functions of $\phi-$mixing sequences, contracting Markov chains, expanding maps of the interval, and symmetric random walks on the circle.

In this paper we are concerned with Cram\'{e}r type moderate deviations  for stationary sequences.
Cram\'{e}r type moderate deviations usually imply the MDP results; see  Fan et al.\ \cite{FGL13} for instance. Furthermore, Cram\'{e}r type moderate deviations  imply  Berry-Esseen bounds; see Corollary \ref{th1.5}. Following the excellent work of Mason and Zhou \cite{MZ12} and Dedecker et al.\ \cite{DMPU09}, we apply our  results to quantile coupling inequalities,  functions of $\phi$-mixing sequences, and contracting Markov chains.

Our approach is based on  martingale approximation and Cram\'{e}r type moderate deviations for martingales due to Fan et al.\ \cite{FGL13}. Cram\'{e}r type moderate deviations for martingales have been established by Ra\v{c}kauskas \cite{Rackauskas90,Rackauskas95},  Grama \cite{G97} and  Grama and Haeusler \cite{GH00,GH06}. Such type of results are very useful for study of stationary sequences, for instance, Wu and Zhao \cite{WZ08}
applied the results of Grama \cite{G97} to establish Cram\'{e}r type moderate deviations for  stationary sequences
with physical dependence measure introduced by Wu \cite{W05}, functionals of linear processes and some nonlinear time series.
See also Cuny and Merlev\`{e}de \cite{CM14} (cf.\ Theorem 3.2 therein) for a result  similar to Wu and Zhao \cite{WZ08}, where Cuny and Merlev\`{e}de \cite{CM14} established a Cram\'{e}r type moderate deviations for  an adapted stationary
sequence in $\mathbf{L}^p$. For relationship among our results and the last two results, we refer to point 3 of Remark \ref{remark2.1}.

The paper is organized as follows. Our main results are stated and discussed in Section \ref{sec2}.
The applications   are given in Section  \ref{sec3}.
Proofs of theorems are deferred to Section \ref{sec4}.

\section{Main results}  \label{sec2}
From now on, assume that the stationary sequence $(X_i)_{i \in \mathbf{Z}}$ is given by $X_i=X_0 \circ T^i, $ where
$T: \Omega \mapsto \Omega$ is a bijective bimeasurable transformation preserving the probability $\mathbf{P}$ on $(\Omega, \mathcal{F})$.
For a subfield $\mathcal{F}_0$ satisfying $\mathcal{F}_0 \subseteq T^{-1}(\mathcal{F}_0)$, let $\mathcal{F}_i= T^{-i}(\mathcal{F}_0).$
Our  theorems  and their corollaries treat  the so-called adapted case, that is $X_0$ being $\mathcal{F}_0$-measurable and so the sequence $(X_i)_{i \in \mathbf{Z}}$ is adapted to the filtration $(\mathcal{F}_i)_{i\in \mathbf{Z}}$. Moreover, we denote the $\mathbf{L}^\infty$-norm by $\|X \|_\infty$, that is the smallest $u$ such that $\mathbf{P}(|X|> u) =0.$

Throughout the paper, let $m=m(n)$ be integers such that $1\leq m \leq   n . $ For instance, we may take
$m=\lfloor n^\alpha \rfloor, \alpha \in (0 , \frac12), $ where $\lfloor x \rfloor$ stands for the largest integer
less than $x.$ Denote
\begin{eqnarray}
&&\varepsilon_m=  \frac m { n^{1/2} \sigma_n } \|X_0\|_\infty, \\
&&\gamma_m=  \frac{1}{m^{1/2} \sigma_n } \sum_{j=1}^{\infty} \frac{1}{ j^{3/2}} \Big \|\mathbf{E}[ S_{mj} | \mathcal{F}_0]\Big \|_\infty  \label{grmman}
\end{eqnarray}
 and
\begin{eqnarray}
\delta_m^2= \frac{1}{m  \sigma_n^2 }  \Big \|\mathbf{E}[ S_m |\mathcal{F}_0]\Big|\!\Big|_\infty^2   +   \Big \| \frac1 {m \sigma_n^2} \mathbf{E}[ S_m^2| \mathcal{F}_0]- 1\Big \|_\infty  ,
\end{eqnarray}
where $\sigma_n =\sqrt{\mathbf{E}  W_n  ^2} >0$.
The following  theorem  gives  a Cram\'{e}r type moderate deviation result for   stationary sequences.

\begin{theorem}\label{th1}
Assume that   $\|X_0\|_\infty < \infty$, and that $X_0$ is $\mathcal{F}_0$-measurable.
Then there exists an absolute constant $\alpha_0 >0 $ such that when  $  \varepsilon_m \leq \frac14, $ $\gamma_m \leq e^{-(80)^2}$ and $  \delta_m^2  +\frac{m}n \leq \alpha_0 $, it holds for  all  $ 0\leq x   \leq \alpha_0  \varepsilon_m^{-1} , $
\begin{eqnarray*}
\Bigg| \ln \frac{\mathbf{P}( W_n  \geq x \sigma_n  )}{1-\Phi \left(  x\right)}  \Bigg|  &\leq&  C_{ \alpha_0 }   \bigg( x^3 \varepsilon_m   + x^2 (  \delta_m^2   +\frac{m}n + \gamma_m|\ln \gamma_m|)  \\
&&\ \ \ \ \ \ \ \  \ \ \ \ \ \ \ + (1+ x  )\big( \varepsilon_m   \left| \ln  \varepsilon_m
 \right|  + \gamma_m|\ln \gamma_m| +    \delta_m    + \sqrt{\frac{m}n}    \big) \bigg)  ,
\end{eqnarray*}
where  $C_{ \alpha_0 }$ depends only on   $\alpha_0.$
In particular, the last inequality implies that
\begin{eqnarray}\label{thls}
 \frac{\mathbf{P}( W_n  \geq x \sigma_n  )}{1-\Phi \left(  x\right)}   = 1 +o(1)
\end{eqnarray}
uniformly  for $\displaystyle 0 \leq x = o(  \min\{ \varepsilon_m^{-1/3}, \, \delta_m^{-1}  \, ,  (n/m)^{1/2}    ,  (\gamma_m |\ln \gamma_m|)^{-1/2} \} )$ as $m \rightarrow \infty.$
Moreover, the same results hold when replacing $ \displaystyle \frac{\mathbf{P}(W_n\geq x \sigma_n)}{1-\Phi \left(  x\right)}$  by $\displaystyle \frac{\mathbf{P}(W_n \leq -x \sigma_n)}{ \Phi \left(-  x\right)}$.
\end{theorem}

\vspace{0.3cm}

\begin{remark}\label{remark2.1}
Let us comment on the results of Theorem  \ref{th1}.
\begin{enumerate}

\item  Assume that
\begin{equation}\label{cond01}
\sum_{n=1}^{\infty} \frac{1}{n^{3/2}}   \Big \|\mathbf{E}[ S_n| \mathcal{F}_0] \Big \|_\infty < \infty,
\end{equation}
and
that there exists  $\sigma  >0 $ such that
\begin{equation}\label{cond02}
 \lim_{n\rightarrow\infty} \Big \| \frac 1 n  \mathbf{E}[ S_n^2| \mathcal{F}_0]- \sigma ^2\Big \|_\infty=0.
\end{equation}
The conditions  (\ref{cond01}) and (\ref{cond02}) were introduced by  Dedecker et al.\ \cite{DMPU09}.
Assume that $m \rightarrow \infty$ and $m/\sqrt{n} \rightarrow 0$ as $n\rightarrow \infty$.
By Lemma 29 of Dedecker et al.\ \cite{DMPU09},   the assumptions of Theorem \ref{th1}   hold  with
$\max\{\varepsilon_m, \gamma_m, \delta_m    \} \rightarrow 0$ as $n \rightarrow \infty.$

\item If $(X_i, \mathcal{F}_i)_{i \in \mathbf{Z}}$ is a   martingale difference sequence,  then
  Theorem \ref{th1} gives a Cram\'{e}r type moderate deviation result  with
  $$ \gamma_m=0   \ \ \  \textrm{and}\ \ \    \displaystyle \delta_m^2=  \Big \| \frac1 {m \sigma_n^2} \sum_{i=1}^m \mathbf{E}[ X_i^2| \mathcal{F}_0]- 1\Big \|_\infty,$$
which is similar to the main theorem of   Grama and  Haeusler \cite{GH00} (see also  Fan et al.\ \cite{FGL13}).
     %


\item  The range of equality  (\ref{thls}) can be very large.  For instance, if\ $\lim _{n\rightarrow \infty} \sigma_n^2=\sigma^2  >0, $ $\big \|\mathbf{E}[ S_n |\mathcal{F}_0]\big \|_\infty =O(1)$ and $\big\| \frac1 n  \mathbf{E}[ S_n^2| \mathcal{F}_0]- \sigma_n^2\big\|_\infty=O\big(\frac1n\big) $
as $n\rightarrow \infty,$ then,  by taking $m=\lfloor n^{2/7} \rfloor,$  equality  (\ref{thls}) holds uniformly  for $ 0 \leq x = o(n^{1/14} /\sqrt{\ln n} )$ as $n \rightarrow \infty.$

\item  For stationary processes,   results similar to Theorem  \ref{th1} can be found in Wu and Zhao \cite{WZ08} and Cuny and Merlev\`{e}de  \cite{CM14}.
Wu and Zhao \cite{WZ08} showed that it is possible to prove the relative error of normal approximation tends to $0$ for
a certain class of stationary processes represented by functions of an i.i.d.\ sequence as soon as the partial sum process can be
well approximated by martingales.  Following the work of Wu and Zhao \cite{WZ08}, Cuny and Merlev\`{e}de (see Theorem 3.2 of \cite{CM14})  proved that under certain conditions for $\mathbf{L}^p$-norm, the relative error of normal approximation tends to $0$ uniformly for $0\leq x = O(\sqrt{\ln n}),$ that is (\ref{thls}) holds  uniformly for $0\leq x = O(\sqrt{\ln n})$. Now Theorem  \ref{th1} shows that the last range could be as large as $0\leq x = o(n^\alpha)$ for some positive constant $\alpha \in (0, \frac12)$ (cf.\ point (iii) of this remark) under  the  conditions for $\mathbf{L}^\infty$-norm (instead of $\mathbf{L}^p$-norm).

\item
The absolute constant $e^{-(80)^2}$ is very small. However, it can be improved to a larger one, provided that
the absolute constant $80$ in the inequality of Peligrad et al.\ \cite{PUW07} (cf.\, inequality (\ref{Peligradcons})) can be improved to a smaller one.

\item  Notice that the quantities $\gamma_m$ and $\delta_m$ can  be estimated via the quantities
  $$\eta_{1, n}:= \sup_{k\geq n} \| \mathbf{E}[X_k |\mathcal{F}_0]\|_\infty\ \ \ \ \ \textrm{and} \ \ \ \ \ \
   \eta_{2, n}:= \sup_{k, l\geq n} \| \mathbf{E}[X_kX_l |\mathcal{F}_0] -\mathbf{E}[X_kX_l ] \|_\infty . $$
Indeed, it is easy to see that
   \begin{eqnarray}
 \gamma_m &\leq&  \frac{1}{m^{1/2} \sigma_n  } \sum_{j=1}^{\infty} \frac{1}{ j^{3/2}}  \Big ( \sum_{i=1}^{mj} \eta_{1, i} \Big )
 \leq \frac{1}{m^{1/2} \sigma_n  }  \sum_{i=1}^{\infty} \eta_{1, i} \sum_{j\geq i/m}  \frac{1}{ j^{3/2}}  \nonumber \\
 &\leq& \frac{C_1}{m^{1/2} \sigma_n  } \Big( \sum_{i=1}^{m} \eta_{1, i}  + \sqrt{m}\sum_{i\geq m}  \frac{\eta_{1, i} }{ i^{1/2}} \Big)  \label{ines35s}
  \end{eqnarray}
  and
   \begin{eqnarray*}
 \delta_m^2 &\leq&  \frac{1}{m  \sigma_n^2  }  \Big[ \Big ( \sum_{i=1}^{m} \eta_{1, i} \Big )^2 +\sum_{i=1}^{m}\| \mathbf{E}[X_i^2 |\mathcal{F}_0] -\mathbf{E}[X_i^2 ] \|_\infty  \\
 && \ \ \ \ \ \ \ \ \ \ \ \ \ \ \ \ \ \  \ \ \ \   +\ 2 \sum_{i=1}^{m-1}\sum_{j=i+1}^{m}\| \mathbf{E}[X_iX_j |\mathcal{F}_0] -\mathbf{E}[X_iX_j ] \|_\infty  \Big],
  \end{eqnarray*}
  where $C_1$ is an absolute constant.
  Splitting the last sum as follows
   \begin{eqnarray*}
  \sum_{1\leq i \leq m/2}    \sum_{  i+1 \leq j \leq 2i }    +  \sum_{1\leq i \leq m/2}    \sum_{  2i+1 \leq j \leq m  }   +   \sum_{ m/2 \leq i \leq m-1}    \sum_{  i+1 \leq j \leq m }  ,
  \end{eqnarray*}
  we infer that
 \begin{eqnarray}
 \delta_m^2 \leq  \frac{C_2}{m  \sigma_n^2  }  \Big[ \Big ( \sum_{i=1}^{m} \eta_{1, i} \Big )^2 +\sum_{1\leq i \leq m/2}  i\eta_{2, i} + \|X_0\|_\infty\sum_{1\leq i \leq m/2}    \sum_{  j \geq 2i  }\eta_{1, j} +   m \!\! \sum_{  i \geq m/2  }  \eta_{2, i}  \Big],  \label{sdc25}
  \end{eqnarray}
  where $C_2$ is an absolute constant.  Moreover,  if $$\lim _{n\rightarrow \infty} \sigma_n^2=\sigma^2  >0 \ \  \ \textrm{and}\  \ \   \max_{i=1,2}\{   \eta_{i, n}  \} =O(  n^{-\beta})$$ for some constant $\beta>1,$  by (\ref{ines35s}) and (\ref{sdc25}),
  then we have  $ \gamma_m  =   O(  m^{-1/2})$   and
  \begin{displaymath}
 \delta_m =  \left\{ \begin{array}{ll}
  O( m^{-1/2} ), & \textrm{\ \ \  if $\beta > 2$,}\\
 O( m^{-1/2} \sqrt{\ln m} ),  & \textrm{\ \ \ if $\beta = 2$,} \\
 O( m^{-(\beta-1)/2} ), & \textrm{\ \ \  if $\beta \in (1, 2)$.}
\end{array} \right.
\end{displaymath}

\item  Assume that $\lim _{n\rightarrow \infty} \sigma_n^2=\sigma^2  >0$. If $  \max_{i=1,2}\{   \eta_{i, n}  \} =O(  n^{-\beta})$ for some constant $\beta\geq 3/2,$ with $m=\lfloor n^{2/7}\rfloor,$ then equality  (\ref{thls}) holds uniformly  for $ 0 \leq x = o(n^{1/14} /\sqrt{\ln n} )$ as $n \rightarrow \infty.$ If $  \max_{i=1,2}\{   \eta_{i, n}  \} =O(  n^{-\beta})$ for some constant $\beta \in (1, 3/2),$ with $m=\lfloor n^{1/(3\beta-1)}\rfloor,$ then equality  (\ref{thls}) holds uniformly  for $ 0 \leq x = o(n^{(\beta-1)/(6\beta-2)}   )$ as $n \rightarrow \infty.$

\end{enumerate}
\end{remark}


 Theorem \ref{th1}   implies the following Berry-Esseen bound.
\begin{corollary}\label{th1.5}
Assume the conditions of Theorem \ref{th1}. Then
\begin{eqnarray}\label{ineq9}
\sup_{ x  } \Big|\mathbf{P}( W_n  \leq  x \sigma_n  )  -  \Phi \left(  x\right)   \Big|   \leq C  \Big(
  \gamma_m|\ln \gamma_m|+ \varepsilon_m \left| \ln  \varepsilon_m
 \right|  + \delta_m  + \sqrt{\frac{m}n}\,  \Big ) ,
\end{eqnarray}
where  $C$ is an absolute constant.
\end{corollary}

\begin{remark}\label{remark2.2}
Let us comment on Corollary \ref{th1.5}.
\begin{enumerate}

\item
 Assume that $\lim _{n\rightarrow \infty} \sigma^2_n = \sigma^2 >0$ and $  \max_{i=1,2}\{   \eta_{i, n}  \}    =O(  n^{-\beta})$ for some constant $\beta>1.$
By point (vi) of Remark \ref{remark2.1},
if $\beta \geq 2,$ then, with $m=\lfloor n^{1/3}\rfloor$, bound (\ref{ineq9}) reaches  its  minimum   of order $n^{-1/6}\ln n .$
If $\beta\in (1, 2),$  then, with $m =\lfloor  n^{1/(\beta+1)} \rfloor$, bound (\ref{ineq9}) gives  its  minimum  of order $n^{-(\beta-1)/(2\beta+2)}\ln n .$

  \item When $(X_i)_{i \in \mathbf{Z}}$ is a uniformly mixing sequence, we refer to Rio \cite{Rio96} for a result similar to Corollary \ref{th1.5}.
   In the paper,  Rio \cite{Rio96} gave a  Berry-Esseen bound of order $n^{-1/2}$ under the condition $ \sum_{k=1}^\infty k \theta_k < \infty,$ where
       $(\theta_k)_{k\geq 1}$ is the sequence of uniformly mixing coefficients.

  \item If $(X_i, \mathcal{F}_i)_{i \in \mathbf{Z}}$ is a stationary  martingale difference sequence,
Corollary \ref{th1.5} gives the following Berry-Esseen bound
\begin{eqnarray}
\sup_{ x  } \Big|\mathbf{P}( W_n  \leq  x \sigma_n  )  -  \Phi \left(  x\right)   \Big|  = O \Big( \frac m { n^{1/2}   } \ln n + \Big \| \frac1 {m \sigma_n^2} \sum_{i=1}^m \mathbf{E}[ X_i^2| \mathcal{F}_0]- 1\Big \|_\infty\Big).
\end{eqnarray}
 When $X_0$ is $\mathbf{L}^p$-bounded (instead of $\mathbf{L}^\infty$-bounded),
Dedecker et al.\,\cite{DMR09} have obtained some rather tight Berry-Esseen bounds.
Notice that Dedecker et al.\,\cite{DMR09} assumed a martingale coboundary decomposition while we do not.
On the other hand Dedecker et al.\,\cite{DMR09} worked in $\mathbf{L}^p$ and we work in $\mathbf{L}^\infty,$
so the results are of independent interest. It is worth noticing that the best rates (for
martingales) provided by Dedecker et al.\,\cite{DMR09} and us are the same.

\end{enumerate}
\end{remark}

Theorem \ref{th1} gives an alternative proof for the following  moderate deviation principle (MDP) result
which is implied by the functional MDP result of  Dedecker et al.\, \cite{DMPU09} under the conditions (\ref{cond01}) and (\ref{cond02}).
\begin{corollary}\label{co0}
Assume the conditions of Theorem \ref{th1}.  Assume   that $\lim _{n\rightarrow \infty} \sigma^2_n = \sigma^2 >0$,  and  that
$\max\{\gamma_m, \delta_m    \} \rightarrow 0$ as $m \rightarrow \infty.$
Let $a_n$ be any sequence of real numbers satisfying $a_n \rightarrow 0$ and $a_n n^{1/2}\rightarrow \infty$
as $n\rightarrow \infty$.  Then  for each Borel set $B \subset \mathbf{R}$,
\begin{eqnarray}
- \inf_{x \in B^o}\frac{x^2}{2 \sigma^2 }  \leq   \liminf_{n\rightarrow \infty}  a_n^2 \ln \mathbf{P}\bigg(a_n W_n      \in B \bigg)
  \leq  \limsup_{n\rightarrow \infty} a_n^2\ln \mathbf{P}\bigg(a_n  W_n     \in B \bigg) \leq  - \inf_{x \in \overline{B}}\frac{x^2}{2  \sigma^2} \, ,   \label{MDP}
\end{eqnarray}
where $B^o$ and $\overline{B}$ denote the interior and the closure of $B$, respectively.
\end{corollary}

 The following  theorem   gives a Bernstein type inequality   for the stationary sequences. Although such type of inequalities are less precise than Cram\'{e}r type moderate deviations, they are available for all positive $x.$ Moreover, they are very useful for establishing quantile coupling inequalities; see Theorem \ref{th3.0}.
\begin{theorem}\label{lemma9}
Assume the conditions of Theorem \ref{th1}. Then  for any $x > 0,$
\begin{eqnarray}
\mathbf{P}\Big(    W_n   \geq x \sigma_n \Big)    &\leq&  \exp\Bigg\{-\frac{(1- \gamma_m|\ln \gamma_m | )^2 x^2}{2 \big(1+\tau_m^2  +\frac{2 }{3 }\varepsilon_m (1- \gamma_m|\ln \gamma_m | ) x \big)}\Bigg\} \, + \, 4 \sqrt{e} \exp  \Bigg\{ - \frac{|\ln \gamma_m |^2 }{2\cdot(81)^2}  \, x^2 \Bigg\},     \label{fsfd03}
\end{eqnarray}
where $\displaystyle \tau_m^2= \delta_m^2 +\frac{m}n   + 4   \varepsilon_m^2.$
\end{theorem}

Assume that $\gamma_m \rightarrow 0$  as $m\rightarrow \infty.$ Then $\gamma_m|\ln \gamma_m | \rightarrow 0$ and $|\ln \gamma_m | \rightarrow \infty $  as $m\rightarrow \infty.$   Thus the second  term  in the r.h.s.\ of (\ref{fsfd03}) is much smaller than the first one for any $x  > 0$ as $m\rightarrow \infty.$
So when $m$ satisfies $m\rightarrow \infty$ and $m/\sqrt{n} \rightarrow 0,$   the bound (\ref{fsfd03}) behaves like $
 \exp\big\{-\frac{  x^2}{2  (1+ \delta_m^2     )}\big\}$ for any $ x> 0 .$

 Next, we apply  Theorems \ref{th1}  and \ref{lemma9} to quantile coupling  inequalities  for  stationary sequences.
We follow  Mason and Zhou \cite{MZ12}, where such type of inequalities have been established
for arbitrary random variables under some Cram\'{e}r type moderate  deviation assumptions.
Using Theorems \ref{th1}, \ref{lemma9}  and Theorem 1 of Mason and Zhou  \cite{MZ12},   we obtain  the following result.
\begin{theorem}\label{th3.0}
Assume the conditions of Theorem \ref{th1}, and that $\gamma_m + \varepsilon_m   + \delta_m +    \sqrt{\frac{m}n} \rightarrow 0$ as $n\rightarrow\infty$.
Let  $\widehat{W}_n= W_n/\sigma_n.$
Then, there exist two positive absolute constants $\alpha$ and $C_{\alpha}$, a standard normal random variable $Z$ and a random variable $Y_n$
can be constructed on a new probability space such that
$Y_n =_d \widehat{W}_n$
and
\begin{eqnarray}\label{fg235}
 |Y_n-Z|  \leq 2  C_\alpha  \Big( Y_n^2\, +1  \Big) \Big (\gamma_m|\ln \gamma_m|+ \varepsilon_m \left| \ln  \varepsilon_m
 \right|  + \delta_m +    \sqrt{\frac{m}n}  \, \Big),
\end{eqnarray}
whenever  \begin{eqnarray}\label{fg2sf}
 |Y_n| \leq  \alpha  \Big(\gamma_m|\ln \gamma_m|+ \varepsilon_m \left| \ln  \varepsilon_m
 \right|  + \delta_m +    \sqrt{\frac{m}n}  \, \Big)^{-1}
\end{eqnarray}
 and $n$ is large enough,  where   $=_d$ stands for equality in distribution.
Furthermore, there exist  two positive absolute constants  $C$ and  $\lambda$  such that for $n$  large enough, we have for all $x\geq0,$
 \begin{eqnarray}\label{fg2sfs5}
 \mathbf{P}\bigg(   \frac{ |Y_n-Z|}{ \gamma_m|\ln \gamma_m| +   \varepsilon_m \left| \ln  \varepsilon_m
 \right|  +  \delta_m + \sqrt{m/n} \,  }  \geq x \bigg) \leq C \exp\Big\{ - \lambda \, x\, \Big\}.
\end{eqnarray}
\end{theorem}

 Assume that $\lim _{n\rightarrow \infty} \sigma^2_n = \sigma^2 >0$ and $  \max_{i=1,2}\{   \eta_{i, n}  \}    =O(  n^{-\beta})$ for some constant $\beta>1.$
By point (i) of Remark \ref{remark2.2},
if $\beta \geq 2,$ then, with $m=\lfloor n^{1/3}\rfloor$, the term $\gamma_m|\ln \gamma_m|+ \varepsilon_m \left| \ln  \varepsilon_m
 \right|  + \delta_m +    \sqrt{\frac{m}n}$ is of order   $n^{-1/6}\ln n .$
If $\beta\in (1, 2),$  then, with $m =\lfloor  n^{1/(\beta+1)} \rfloor$, the term $\gamma_m|\ln \gamma_m|+ \varepsilon_m \left| \ln  \varepsilon_m
 \right|  + \delta_m +    \sqrt{\frac{m}n}$ is of order $n^{-(\beta-1)/(2\beta+2)}\ln n .$

\section{Applications}\label{sec3}
In this section, we present some applications of our results. For more interesting applications, such as  expanding map and symmetric random walk on the circle, we refer to Corollary 18 and Proposition 20 of Dedecker et al.\ \cite{DMPU09}. Under their corresponding conditions, the conditions of Theorem \ref{th1} hold.

\subsection{$\phi$-mixing sequences}
Let $Y$ be a random variable with values in a Polish space $\mathcal{Y}.$ If $\mathcal{M}$
is a $\sigma$-field,  the $\phi$-mixing coefficient between $\mathcal{M}$ and $\sigma(Y)$ is defined by
\begin{equation}\label{phi}
\phi(\mathcal{M}, \sigma(Y))=\sup_{A \in \mathfrak{B}(\mathcal{Y})} \Big\| \mathbf{P}_{Y|\mathcal{M}}(A) -\mathbf{P}_Y(A) \Big\|_\infty.
\end{equation}
For  a sequence of random variables $(X_i)_{i \in \mathbf{Z}}$ and a positive integer $m,$ denote
\[
\phi_m(n)=\sup_{i_m> ...> i_1\geq n}\phi(\mathcal{F}_0, \sigma(X_{i_1},...,X_{i_m} )),
 \]
 and let $\phi(k)=\lim_{m\rightarrow \infty}\phi_m(k)$ be the usual $\phi$-mixing coefficient.
Under the following conditions
\begin{equation}\label{phic01}
 \sum_{k\geq 1}k^{1/2}\phi_1(k) < \infty  \ \   \ \ \textrm{and }\ \ \ \ \lim_{k \rightarrow \infty} \phi_2(k)=0 ,
\end{equation}
 Dedecker et al.\ \cite{DMPU09}   obtained a MDP result for bounded random variables.
See also Gao \cite{G96} for an earlier MDP result under a condition stronger than  (\ref{phic01}), that is $ \sum_{k\geq1} \phi(k) < \infty$.

  When the random variables $(X_i)_{i \in \mathbf{Z}}$ are  bounded, it holds $\eta_{1, n}=O(\phi_1(n))$ and $\eta_{2, n}=O(\phi_2(n))$ as $n\rightarrow \infty$. By point (vii) of Remark \ref{remark2.1},  we have the following result.
 \begin{proposition}   Assume that the random variables $(X_i)_{i \in \mathbf{Z}}$ are  bounded, $\lim _{n\rightarrow \infty} \sigma_n^2= \sigma^2 >0$ and
 \[
    \max_{i=1,2}\{   \phi_i(n)  \} =O ( n^{-\beta}),\ \    \ \  n\rightarrow \infty,
 \]
 for some constant $\beta>1.$
 \begin{description}
  \item[\textbf{[i]}]  If $\beta \geq 3/2$, then    (\ref{thls}) holds uniformly  for $ 0 \leq x = o(n^{1/14} /\sqrt{\ln n} )$ as $n \rightarrow \infty.$
  \item[\textbf{[ii]}]   If $\beta\in (1, 3/2)$, then   (\ref{thls}) holds uniformly  for $ 0 \leq x = o(n^{(\beta-1)/(6\beta-2)}   )$ as $n \rightarrow \infty.$
\end{description}
\end{proposition}

\subsection{Functions of $\phi$-mixing sequences}\label{ap3.2}
Let $(\varepsilon_i)_{i \in \mathbf{Z}}=(\varepsilon_0 \circ T^i)_{i \in \mathbf{Z}}$ be a stationary sequence of $\phi$-mixing random variables
taking values in a subset $A$ of a Polish space $ \mathcal{X}.$ Denote by $\phi_\varepsilon(n)$ the coefficient
\[
\phi_\varepsilon(n)= \phi( \sigma(\varepsilon_i, i\leq 0), \sigma(\varepsilon_i, i\geq n)  ),
\]
 where  $\phi$ is defined by (\ref{phi}).
 Let $H$ be a function from $A^{\mathbf{N}}$ to $\mathbf{R}$ satisfying the following condition
\begin{description}
  \item[(A):] \ \   $ \textrm{for any}\  i \geq 0,    \  \ \  \ \sup_{x \in A^{\mathbf{N}}, \, y \in A^{\mathbf{N}} }\Big| H(x)-H(x^{(i)}y) \Big|  \leq R_i,\ \ \textrm{where}\ \ R_i \ \textrm{decreases to }\ 0,$
\end{description}
where the sequence  $x^{(i)}y$ is defined by $(x^{(i)}y)_j=x_j$ for $j<i$ and $(x^{(i)}y)_j=y_j$ for $j\geq i.$ Define the stationary sequence $X_k=X_0\circ T^k$
by
\begin{equation}\label{thedfx}
X_k=H((\varepsilon_{k-i})_{i \in \mathbf{N}} ) -\mathbf{E}[H((\varepsilon_{k-i})_{i \in \mathbf{N}})].
\end{equation}

 Dedecker et al.\ \cite{DMPU09} gave a MDP result for $(X_k)_{k\geq1}$, see Propositions 12 therein.
 From the proof of Propositions 12 of  \cite{DMPU09}, it is easy to see that
 $$ \max_{i=1,2}\{\eta_{i, n}  \} =O \Big(   R_n + \sum_{i=1}^n R_{n-i} \phi_\varepsilon(i)  \Big)   .$$
 Notice that when $\sigma^2:= \sum_{k \in \mathbf{Z}}\mathbf{E}[X_0X_k] >0,$ it holds
$ \lim_{n\rightarrow \infty} \sigma_n^2 = \sigma^2.$
By point (vii) of Remark \ref{remark2.1}, we have the following Cram\'{e}r type  moderate deviations.
\begin{proposition}\label{Pro3.1s}
  Let $(X_k)_{k \in \mathbf{Z}}$ be defined by (\ref{thedfx}), for a function $H$ satisfying condition (A). Assume
 \begin{equation}\label{sfsf02}
  R_n + \sum_{i=1}^n R_{n-i} \phi_\varepsilon(i) = O(  n^{-\beta}) , \ \ \ \ n\rightarrow\infty,
 \end{equation}
 for some constant $\beta>1,$  and   $\sigma^2:= \sum_{k \in \mathbf{Z}}\mathbf{E}[X_0X_k] >0.$
 \begin{description}
  \item[\textbf{[i]}]  If $\beta \geq 3/2$, then     (\ref{thls}) holds uniformly  for $ 0 \leq x = o(n^{1/14} /\sqrt{\ln n} )$ as $n \rightarrow \infty.$
  \item[\textbf{[ii]}]   If $\beta \in (1, 3/2)$, then   (\ref{thls}) holds uniformly  for $ 0 \leq x = o(n^{(\beta-1)/(6\beta-2)}   )$ as $n \rightarrow \infty.$
\end{description}
\end{proposition}

\subsection{Contracting Markov chains}\label{ap3.3}
Let $(Y_n)_{n\geq 0} $ be a stationary Markov chain of bounded random variables with invariant measure $\mu$ and transition kernel $K.$
Denote by $\| \cdot \|_{\infty, \mu} $ the essential norm with respect to $\mu.$ Let $\Lambda_1$ be the set of $1$-Lipschitz
functions. Assume that the chain satisfies the  following condition:
\begin{description}
  \item[(B):] \ \ there exist two constants $C>0$ and $\rho \in (0, 1)$ such that \[\sup_{g \in \Lambda_1} \|  K^n(g) - \mu(g)  \|_{\infty, \mu} \leq C \rho^n,   \]
    \quad \quad \  and for any $m\geq 0,$
  \[ \sup_{f, g \in \Lambda_1}  \Big\|  K^n\big(fK^m(g) \big) - \mu\big(fK^m(g)\big)  \Big\|_{\infty, \mu}   \leq C \rho^n. \]
\end{description}
We shall see in the next proposition that MDP result holds  for the sequence
 \begin{equation} \label{fds}
 X_n=f(Y_n) - \mu(f)
 \end{equation}
as soon as the function $f$ belongs to the class $\mathcal{L}$ introduced by Dedecker et al.\ \cite{DMPU09}. Let $\mathcal{L}$ be the class of functions
$f:  \mathbf{R} \mapsto \mathbf{R}$ such that $|f(x) -f(y)| \leq g(|x-y|)$, where $g$ is a  concave and non-decreasing function and
satisfies
\begin{equation} \label{fsddsdfs}
 \int_0^1 \frac{g(t)}{ t \sqrt{|\ln t| }}  dt < \infty.
 \end{equation}
Clearly,  (\ref{fsddsdfs}) holds if $g(t) \leq c |\ln(t)|^{-\gamma}$ for some constants $c>0$ and  $\gamma > 1/2.$
In particular,  $\mathcal{L}$ contains the class of $\alpha$-H\"{o}lder continuous  functions from $[0, 1]$ to $\mathbf{R}$, where $\alpha \in (0, 1].$

Dedecker et al.\ \cite{DMPU09} gave a MDP result for $(Y_n)_{n\geq0}$, see Propositions 14 therein.
From the proof of Propositions 14  of  \cite{DMPU09}, it is easy to see that
 $$ \max_{i=1,2}\{\eta_{i, n}  \} =O \big(  g(C \rho^n ) \big)  ,$$
 where $C$ is given by condition (B).
\begin{proposition}\label{ap02s}
Assume that the stationary Markov chain $(Y_n)_{n\geq 0}$ satisfies condition  (B), and let $X_n$
be defined by (\ref{fds}).
Assume $f \in \mathcal{L},$
\[
\sigma^2:=\sigma^2(f)=\mu \Big((f-\mu(f))^2 \Big) + 2 \sum_{n>0} \mu \Big(K^n(f)\cdot (f-\mu(f))\Big)  >0
\]
and
 \begin{equation} \label{dfs456}
 g(C \rho^n ) = O(  n^{-\beta}) , \ \ \ \ n\rightarrow\infty,
 \end{equation}
 for some constant $\beta>1.$
\begin{description}
  \item[\textbf{[i]}]  If $\beta \geq 3/2$, then  (\ref{thls}) holds uniformly  for $ 0 \leq x = o(n^{1/14} /\sqrt{\ln n} )$ as $n \rightarrow \infty.$
  \item[\textbf{[ii]}]   If $\beta \in (1, 3/2)$, then  (\ref{thls}) holds uniformly  for $ 0 \leq x = o(n^{(\beta-1)/(6\beta-2)}   )$ as $n \rightarrow \infty.$
\end{description}
\end{proposition}

Notice that if $g(t) \leq D |\ln(t)|^{-\beta}$ for some constants $D>0$ and  $\beta > 1,$ then   (\ref{dfs456}) is satisfied.

%

\section{Proofs of Theorems and Corollaries }\label{sec4}
The proofs of our results are mainly based on the following lemmas, which give some exponential deviation inequalities
for the partial sums of dependent random variables.

\subsection{Preliminary lemmas}\label{sec4.0}
Let $(\xi_i,\mathcal{F}_i)_{i=0,...,n}$ be a sequence of martingale differences, defined on some
 probability space $(\Omega ,\mathcal{F},\mathbf{P})$,  where $\xi_0=0 $,  $\{\emptyset, \Omega\}=\mathcal{F}_0\subseteq ...\subseteq \mathcal{F}_n\subseteq
\mathcal{F}$ are increasing $\sigma$-fields. Set
\begin{equation}
M_{0}=0,\ \ \ \ \ M_k=\sum_{i=1}^k \xi_i,\quad k=1,...,n.  \label{xk}
\end{equation}
Then $M=(M_k,\mathcal{F}_k)_{k=0,...,n}$ is a martingale.  Denote
  $\left\langle M\right\rangle $  the quadratic characteristic of the
martingale $M$, that is
\begin{equation}\label{quad}
\left\langle M\right\rangle _0=0,\ \ \ \ \ \left\langle M\right\rangle _k=\sum_{i=1}^k\mathbf{E}[\xi_i^2|\mathcal{F}
_{i-1}],\quad k=1,...,n.
\end{equation}

 Assume   the following two conditions:
\begin{description}
\item[(C1)]  There exists   $\epsilon_n \in (0, \frac12]$ such that
\[
\Big|\mathbf
{E}[\xi_{i}^{k}  | \mathcal{F}_{i-1}] \Big| \leq \frac12 k!\epsilon_n^{k-2} \mathbf{E}[\xi_{i}^2 | \mathcal{F}_{i-1}],\ \ \ \ \ \textrm{for all}\ k\geq 3\, \ \textrm{and}\ \, 1\leq i\leq n;
\]
\item[(C2)]  There exists   $ \iota_n \in [0, \frac12]$ such that
$ \| \left\langle M\right\rangle _n-1\|_\infty \leq  \iota_n^2. $
\end{description}
Clearly, condition (C1) is satisfied for bounded martingale differences $\|\xi_{i}\|_\infty \leq \epsilon_n.$

In the proof  of  Theorem \ref{th1}, we need the following Cram\'{e}r moderate deviation expansions for martingales, which is a simple consequence of Theorems 2.1 and 2.2  of Fan et al.\ \cite{FGL13}.
\begin{lemma}\label{lemma4.1}
Assume conditions (C1) and (C2). Then there is an absolute constant $\alpha_0  >0$ such that
for all $0\leq x \leq \alpha_0  \,\epsilon_n^{-1}$ and $\iota_n \leq \alpha_0$,
\begin{equation} \label{t1ie1}
\bigg| \ln \frac{\mathbf{P}(M_n \geq x)}{1-\Phi \left( x\right)} \bigg| \leq C_{\alpha_0 } \! \Big( x^3 \epsilon_n  + x^2 \iota_n^2 + (1+ x) \left( \epsilon_n \left| \ln  \epsilon_n
 \right|+  \iota_n \right)\Big)  ,
\end{equation}
where $C_{\alpha_0 }$ depends only on $\alpha_0 .$   Moreover, the same equality remains true when $\frac{\mathbf{P}(M_n\geq x)}{1-\Phi \left(  x\right)}$ is replaced by $\frac{\mathbf{P}(M_n \leq -x)}{ \Phi \left(-  x\right)}$.
\end{lemma}

In the proof of Theorem \ref{lemma9}, we make use of the following Freedman inequality \cite{FR75}. 
\begin{lemma}\label{lemma4.2}
 Assume that $\xi_{i} \leq a$ for some constant $a$ and all $1\leq i \leq n.$   Then
for all $ x \geq 0$ and $v_n > 0$,
\begin{eqnarray}
  \mathbf{P}\left( M_n \geq x\ \mbox{and}\ \langle M \rangle_{n}\leq v_n^2\   \right) \leq \exp\bigg\{-\frac{x^2}{2(v_n^2+\frac{a}{3}x )}\bigg\}.\label{fgl1}
\end{eqnarray}
\end{lemma}

We also use  the following  exponential inequality of  Peligrad et al.\ \cite{PUW07} (cf.\ Proposition 2 therein), which plays an important role in the proof  of  Theorem \ref{lemma9}.
\begin{lemma}\label{lemma1}
  Let $(X_i)_{i \in \mathbf{Z}}$ be a stationary sequence of random variables adapted to the filtration $(\mathcal{F}_i)_{i \in \mathbf{Z}}$. Then for all $x \geq 0,$
  \begin{equation}\label{Peligradcons}
\mathbf{P}\bigg( \max_{1\leq i \leq n} |S_i| \geq x \bigg)  \leq 4 \sqrt{e} \exp  \bigg\{ - \frac{x^2 }{ 2n (  \|X_1\|_\infty + 80 \sum_{j=1}^n j^{-3/2}\|\mathbf{E}[S_j|\mathcal{F}_0]\|_\infty )^2 } \bigg\}.
\end{equation}
\end{lemma}
\subsection{Proof of Theorem  \ref{th1}}\label{sec4.1}
Let   $k=k(n,m)=\lfloor n/m \rfloor$ be the integer part of $n/m$. The initial step of the proof is to divide the random variables into blocks of size $m$ and
to make the sums in each block
\[
X_{i,m}=\sum_{j=(i-1)m+1}^{im} X_j,\ \ 1\leq i  \leq k,\ \ \ \textrm{and} \ \ \ X_{k+1,m}=\sum_{j= km+1}^{n} X_j.
\]
It is easy to see that $S_n= \sum_{i=1}^{k+1} X_{i,m}. $
Define
\[
 D_{i,m}=X_{i,m} - \mathbf{E}[X_{i,m} | \mathcal{F}_{(i-1)m} ],\ \  1\leq i  \leq k.
 \]
Then $(D_{i,m}, \mathcal{F}_{im} )_{1\leq i \leq k}$ is a stationary  sequence of bounded    martingale differences, that is
\[
\|D_{i,m}\|_\infty \leq 2m \|X_0 \|_\infty.
\]
Notice that
\[
\mathbf{E}[ D_{i,m}^2 | \mathcal{F}_{(i-1)m}  ]=\mathbf{E}[ X_{i,m}^2 | \mathcal{F}_{(i-1)m}  ] -(\mathbf{E}[ X_{i,m}  | \mathcal{F}_{(i-1)m}  ])^2,
\]
and that, by stationarity, it follows that
\[
 \frac1n \Big\| \sum_{i=1}^k (\mathbf{E}[ X_{i,m}  | \mathcal{F}_{(i-1)m}  ])^2 \Big\|_\infty \leq   \frac1m   \Big\| \mathbf{E}[ S_{m}  | \mathcal{F}_{0}  ]\Big\|_\infty^2 .
 \]
Moreover,
\begin{eqnarray*}
\Big\| \frac1n \sum_{i=1}^k  \mathbf{E}[ X_{i,m}^2  | \mathcal{F}_{(i-1)m}  ] -   \sigma_n^2 \Big\|_\infty
&\leq&   \frac1n \sum_{i=1}^k  \Big\| \mathbf{E}[ X_{i,m}^2 | \mathcal{F}_{(i-1)m}  ] -  m \sigma_n^2\Big\|_\infty^2 + \frac{n-mk}n \sigma_n^2  \\
 &\leq&  \Big\|\frac1m \mathbf{E}[ S_m ^2 | \mathcal{F}_{0}  ] - \sigma_n^2\Big\|_\infty^2 +  \frac{m}n \sigma_n^2 .
\end{eqnarray*}
Consequently, it holds
\begin{eqnarray*}
&& \Big\| \frac1n \sum_{i=1}^{k }\mathbf{E}[ D_{i,m}^2 | \mathcal{F}_{(i-1)m}  ] - \sigma_n^2 \Big\|_\infty \\
&&\leq  \Big\| \frac1n \sum_{i=1}^k \mathbf{E}[ X_{i,m}^2 | \mathcal{F}_{(i-1)m}  ] -   \sigma_n^2\Big\|_\infty^2      \ +\  \frac1n \Big\| \sum_{i=1}^k (\mathbf{E}[ X_{i,m}  | \mathcal{F}_{(i-1)m}  ])^2 \Big\|_\infty   \\
 && \leq \Big\|\frac1m \mathbf{E}[ S_m^2  | \mathcal{F}_{0}  ] - \sigma_n^2\Big\|_\infty^2 +     \frac{m}n \sigma_n^2  \ +\   \frac1m   \Big\| \mathbf{E}[ S_{m}  | \mathcal{F}_{0}  ]\Big\|_\infty^2 \\
 && = \big( \delta_m^2 +\frac{m}n \big) \sigma_n^2
\end{eqnarray*}
and
\[
\|n^{-1/2} D_{i,m}\|_\infty \leq 2  \sigma_n \varepsilon_m  .
\]
Denote  $\xi_i =D_{i,m} /( n^{1/2}\sigma_n)  $ and $ M_k=\sum_{i=1}^k \xi_i .$ Then it is obvious that  $$|\xi_i| \leq 2 \varepsilon_m    \ \ \ \ \textrm{and}\ \ \ \ \| \langle M\rangle_k  -1 \|_{\infty} \leq    \delta_m^2  +\frac{m}n .$$
Assume $   \varepsilon_m \leq \frac14 $ and $  \delta_m^2  +\frac{m}n \leq   \alpha_0$,
where $\alpha_0 \in (0, \frac12]$ is given by Lemma \ref{lemma4.1}.
By Lemma \ref{lemma4.1}, we have for all $0\leq x  \leq \alpha_0   \varepsilon_m^{-1}$,
\begin{eqnarray}
   \bigg|\ln \frac{\mathbf{P}( M_k \geq x )}{1-\Phi \left(  x\right)}  \bigg|  \leq  C'_{\alpha_0 }  \, \bigg( x^3 \varepsilon_m     + x^2 (   \delta_m^2  +\frac{m}n)  + (1+ x)\big(  \varepsilon_m  \left| \ln   \varepsilon_m
 \right|  +    \delta_m       + \sqrt{ m/n}  \,  \big)\bigg)     ,\label{ineq4}
\end{eqnarray}
where $C'_{\alpha_0}$  depends only on  $\alpha_0.$
Notice that for all $x\geq 0$ and $|\varepsilon|  \leq \frac 1 2$,
\begin{eqnarray}
\frac{1-\Phi \left( x + \varepsilon \right)}{1-\Phi \left( x\right) }=   \exp\Big\{ \theta  \sqrt{2 \pi} (1+x) |\varepsilon| \Big\} \label{hjklms}
\end{eqnarray}
and
\begin{eqnarray*}
\frac1{\sqrt{n}   \sigma_n }\big\|X_{k+1,m}  \big\|_{\infty} \leq  \frac1{\sqrt{n}   \sigma_n }   (n- km)  \big\|X_{0}  \big\|_{\infty}   \leq      \varepsilon_m ,
\end{eqnarray*}
where $|\theta|\leq 1.$
It is obvious that
\[
M_k +   \frac1{\sqrt{n}   \sigma_n } X_{k+1,m}     = \frac1{\sqrt{n}  \sigma_n  }  \Big( S_n - \sum_{i=1}^{k } \mathbf{E}[X_{i,m} | \mathcal{F}_{(i-1)m} ] \Big) .
\]
Therefore, by (\ref{ineq4}) and (\ref{hjklms}),  for all  $0\leq x  \leq \alpha_0   \varepsilon_m^{-1} ,$
\begin{eqnarray*}
 &&\frac{\mathbf{P}( S_n - \sum_{i=1}^{k } \mathbf{E}[X_{i,m} | \mathcal{F}_{(i-1)m} ] \geq x\sigma_n n^{1/2})}{1-\Phi \left(  x\right)} \nonumber \\
 &&\  \leq \  \frac{\mathbf{P}( M_k   \geq x +     \varepsilon_m     )}{1-\Phi \left(  x  +     \varepsilon_m   \right)} \cdot
 \frac{ 1-\Phi \left(  x  +     \varepsilon_m    \right)  }{ 1-\Phi \left(  x \right) } \nonumber \\
 && \  \leq \ \exp\bigg\{   C _{\alpha_0 }  \, \bigg( x^3 \varepsilon_m      + x^2 ( \delta_m^2       +\frac{m}n)  + (1+ x  )\big( \varepsilon_m  \left| \ln  \varepsilon_m
 \right|  +   \delta_m        + \sqrt{m/n}    \big)\bigg)  \bigg\} .
\end{eqnarray*}
Similarly, we have for all  $0\leq x  \leq \alpha_0   \varepsilon_m^{-1} ,$
\begin{eqnarray*}
 &&\frac{\mathbf{P}( S_n - \sum_{i=1}^{k } \mathbf{E}[X_{i,m} | \mathcal{F}_{(i-1)m} ] \geq x\sigma_n n^{1/2})}{1-\Phi \left(  x\right)} \nonumber \\
 && \  \geq \ \exp\bigg\{ -  C _{\alpha_0 }  \, \bigg( x^3 \varepsilon_m     + x^2 (  \delta_m^2    +\frac{m}n)  + (1+ x  )\big( \varepsilon_m \left| \ln  \varepsilon_m
 \right|  +    \delta_m        + \sqrt{ m/n}    \big)\bigg)  \bigg\} .
\end{eqnarray*}
The last two inequalities imply that  for all  $0\leq x  \leq \alpha_0   \varepsilon_m^{-1} ,$
\begin{eqnarray}
 && \bigg|\ln \frac{\mathbf{P}( S_n - \sum_{i=1}^{k } \mathbf{E}[X_{i,m} | \mathcal{F}_{(i-1)m} ] \geq x\sigma_n n^{1/2})}{1-\Phi \left(  x\right)} \bigg| \nonumber \\
 && \  \leq \   C _{\alpha_0 }  \, \bigg( x^3 \varepsilon_m   + x^2 (  \delta_m^2   +\frac{m}n)  + (1+ x  )\big( \varepsilon_m   \left| \ln  \varepsilon_m
 \right|  +    \delta_m        + \sqrt{ m/n}    \big)\bigg)    .\label{ineq5}
\end{eqnarray}
By Lemma \ref{lemma1},  we derive that for all $x \geq 0$,
\begin{eqnarray}
&& \mathbf{P} \Big(  \Big|  \sum_{i=1}^{k } \mathbf{E}[X_{i,m} | \mathcal{F}_{(i-1)m} ] \Big| \geq  x \sigma_n n^{1/2} \Big) \nonumber  \\
 &&\leq
 4 \sqrt{e} \exp  \bigg\{ - \frac{n \sigma_n^2 x^2 }{ 2 k ( \big\| \mathbf{E}[ S_{m}  | \mathcal{F}_{0}  ]\big\|_\infty + 80 \sum_{j=1}^k j^{-3/2}\|\mathbf{E}[S_{jm}|\mathcal{F}_0]\|_\infty )^2 } \bigg\}  \nonumber\\
  &&\leq
 4 \sqrt{e} \exp  \bigg\{ - \frac{x^2   }{ 2 \cdot(81)^2   \gamma_m^2   } \bigg\}.  \label{ineq6}
\end{eqnarray}
It is easy to see that for all $x \geq 0$,
\begin{eqnarray}
 \mathbf{P}\Big( W_n   \geq x\sigma_n  \Big)& \leq &   \mathbf{P}\bigg( S_n - \sum_{i=1}^{k } \mathbf{E}[X_{i,m} | \mathcal{F}_{(i-1)m} ] \geq (1-\gamma_m|\ln \gamma_m|)x\sigma_n n^{1/2} \bigg)\nonumber \\
 & & \  + \ \mathbf{P}\bigg(  \sum_{i=1}^{k } \mathbf{E}[X_{i,m} | \mathcal{F}_{(i-1)m} ] \geq  \gamma_m|\ln \gamma_m| x\sigma_n n^{1/2} \bigg).\label{ineq7}
\end{eqnarray}
By the inequalities (\ref{ineq5})-(\ref{ineq7}), it follows that for  all $0\leq x  \leq \alpha_0 \varepsilon_m^{-1} ,$
\begin{eqnarray*}
\frac{\mathbf{P}( W_n   \geq x\sigma_n)}{1-\Phi \left(  x\right)}  &\leq &\frac{1-\Phi \left( (1-\gamma_m|\ln \gamma_m|)x\right)}{1-\Phi \left(  x\right) } \\
&& \times \exp\bigg\{   C _{\alpha_0 }  \, \bigg( x^3 \varepsilon_m   + x^2 (  \delta_m^2   +\frac{m}n)  + (1+ x  )\big( \varepsilon_m   \left| \ln  \varepsilon_m
 \right|  +    \delta_m        + \sqrt{\frac{m}n}    \big)\bigg)  \bigg\} \\
 && \ \   +\ \frac{4 \sqrt{e} }{1-\Phi \left(  x\right) }\exp  \bigg\{ - \frac{ 1}{ 2 \cdot(81)^2   } (\ln \gamma_m )^2 x^2 \bigg\}.
\end{eqnarray*}
Using the following two-sided bound  on tail probabilities of the standard normal random variable
\begin{eqnarray}\label{fgsgj1}
\frac{1}{\sqrt{2 \pi}(1+x)} e^{-x^2/2} \leq 1-\Phi ( x ) \leq \frac{1}{\sqrt{ \pi}(1+x)} e^{-x^2/2}, \ \   x\geq 0,
\end{eqnarray}
we deduce that  for   all $\gamma_m \leq e^{-(80)^2}$ and $ 1 \leq x   \leq \alpha_0 \varepsilon_m^{-1}, $
\begin{eqnarray}
 \frac{\mathbf{P}(W_n   \geq x\sigma_n)}{1-\Phi \left(  x\right)}   &\leq& \exp\bigg\{   C_{ \alpha_0}  \bigg( x^3 \varepsilon_m   + x^2 (  \delta_m^2   +\frac{m}n + \gamma_m|\ln \gamma_m|)  + (1+ x  )\big( \varepsilon_m   \left| \ln  \varepsilon_m
 \right|  +    \delta_m        + \sqrt{\frac{m}n}    \big)\bigg)   \bigg\}  \nonumber \\
 && \ + \ C_1    \times  \exp  \bigg\{ - \frac{ 1}{ 4 \cdot(81)^2   } |\ln \gamma_m | x^2 \bigg\} \nonumber \\
 &\leq& \exp\bigg\{   C_{ \alpha_0}   \bigg( x^3 \varepsilon_m   + x^2 (  \delta_m^2   +\frac{m}n + \gamma_m|\ln \gamma_m|)  + (1+ x  )\big( \varepsilon_m   \left| \ln  \varepsilon_m
 \right|  +    \delta_m        + \sqrt{\frac{m}n}    \big)\bigg)   \bigg\}  \nonumber \\
 && \ + \ C_2 \gamma_m |\ln \gamma_m |  x^2     \nonumber \\
   &\leq&  \exp\bigg\{  C'_{  \alpha_0}   \Big(x^3 \varepsilon_m   + x^2 (  \delta_m^2   +\frac{m}n + \gamma_m|\ln \gamma_m|)  + (1+ x  )\big( \varepsilon_m   \left| \ln  \varepsilon_m
 \right|  +    \delta_m        + \sqrt{\frac{m}n}    \big)   \Big)  \bigg\} .  \nonumber  \\
 &&  \label{need01}
\end{eqnarray}
Notice that for $x\geq 0,$
\begin{eqnarray}
 \mathbf{P}\Big( W_n   \geq x\sigma_n  \Big)   & \geq &     \mathbf{P}\bigg( S_n - \sum_{i=1}^{k } \mathbf{E}[X_{i,m} | \mathcal{F}_{(i-1)m} ] \geq (1+ \gamma_m|\ln \gamma_m|)x\sigma_n n^{1/2} \bigg)\nonumber \\
 & & - \ \mathbf{P}\bigg(  \sum_{i=1}^{k } \mathbf{E}[X_{i,m} | \mathcal{F}_{(i-1)m} ] \leq -   \gamma_m|\ln \gamma_m| x\sigma_n n^{1/2} \bigg).
\end{eqnarray}
By an argument similar to the proof  of (\ref{need01}),   we  deduce that  for  all  $ 1 \leq x   \leq \alpha_0 \varepsilon_m^{-1} ,$
\begin{eqnarray}
&& \frac{\mathbf{P}( W_n   \geq x\sigma_n)}{1-\Phi \left(  x\right)} \nonumber \\
  &&\ \ \ \ \   \geq   \exp\bigg\{ - C'_{  \alpha_0}   \bigg(x^3 \varepsilon_m   + x^2 (  \delta_m^2   +\frac{m}n + \gamma_m|\ln \gamma_m|)  + (1+ x  )\big( \varepsilon_m   \left| \ln  \varepsilon_m
 \right|  +    \delta_m        + \sqrt{\frac{m}n}    \big)  \bigg)  \bigg\}    . \label{need02}
\end{eqnarray}
Combining (\ref{need01}) and (\ref{need02}) together, we obtain the desired equality for  all  $ 1 \leq x   \leq \alpha_0 \varepsilon_m^{-1} .$
Next, we consider the case where $ x \in [0,  1].$
Notice that (\ref{ineq5})  holds also for $(-X_i)_{i \in \mathbf{Z}}.$ Thus, from (\ref{ineq5}), we have
 \begin{eqnarray}
   \sup_{|x| \leq 2} \Big | \mathbf{P}\Big( S_n - \sum_{i=1}^{k } \mathbf{E}[X_{i,m} | \mathcal{F}_{(i-1)m} ] \geq x\sigma_n n^{1/2} \Big)   - \big(  1-\Phi \left(x\right) \big)  \Big|     \leq    C_{\alpha_0} \, \Big(   \varepsilon_m  | \ln  \varepsilon_m |   + \delta_m + \sqrt{ m/n}  \Big)   . \ \ \ \ \
\end{eqnarray}
 For all  $  x \in [0,  1], $  we deduce that
\begin{eqnarray*}
&&  \mathbf{P}\Big( W_n   \geq x\sigma_n  \Big) - \Big( 1- \Phi \left(  x\right) \Big)     \\
 && \geq  \mathbf{P}\Big( S_n - \sum_{i=1}^{k } \mathbf{E}[X_{i,m} | \mathcal{F}_{(i-1)m} ] \geq  ( x   - \gamma_m|\ln \gamma_m| ) \sigma_n n^{1/2}    \Big)   -  \Big( 1- \Phi \left(  x\right) \Big)   \\
 && \  \ \ \ - \ \mathbf{P}\bigg(   \sum_{i=1}^{k } \mathbf{E}[X_{i,m} | \mathcal{F}_{(i-1)m} ] \geq  \gamma_m|\ln \gamma_m|  \sigma_n n^{1/2} \bigg)\\
&& \geq - C_{\alpha_0} \, \Big(   \varepsilon_m  | \ln  \varepsilon_m |   + \delta_m + \sqrt{\frac{m}n} \Big) - \Big|   \Big( 1- \Phi \left(  x  - \gamma_m|\ln \gamma_m|\right) \Big)  -  \Big( 1- \Phi \left(  x\right) \Big)  \Big| \\
 && \  \ \ \ - \ \mathbf{P}\bigg(  \sum_{i=1}^{k } \mathbf{E}[X_{i,m} | \mathcal{F}_{(i-1)m} ] \geq  \gamma_m|\ln \gamma_m|  \sigma_n n^{1/2} \bigg)\\
 && \geq  -C_{\alpha_0} \, \Big(   \varepsilon_m  | \ln  \varepsilon_m | + \gamma_m|\ln \gamma_m| + \delta_m + \sqrt{ m/n}  \Big),
\end{eqnarray*}
where the last line follows by (\ref{ineq6}).
Similarly, we have for all  $  x \in [0,  1], $
\begin{eqnarray*}
 \mathbf{P}\Big( W_n   \geq x\sigma_n  \Big) - \Big( 1- \Phi \left(  x\right) \Big)  \leq   C_{\alpha_0} \, \Big(   \varepsilon_m  | \ln  \varepsilon_m | + \gamma_m|\ln \gamma_m| + \delta_m + \sqrt{ m/n}  \Big).
\end{eqnarray*}
The last two inequalities imply that for all  $  x \in [0,  1], $
\begin{eqnarray*}
  \Big|\mathbf{P}\Big( W_n   \geq x\sigma_n  \Big) - \Big( 1- \Phi \left(  x\right) \Big)  \Big|  \leq  C_{\alpha_0} \, \Big(   \varepsilon_m  | \ln  \varepsilon_m | + \gamma_m|\ln \gamma_m| + \delta_m + \sqrt{ m/n}  \Big).
\end{eqnarray*}
The last inequality implies   the desired equality for  all  $ x \in [0,  1].$

Since $(-X_i)_{i \in \mathbf{Z}}$ also satisfies the conditions of Theorem  \ref{th1},
the same equalities remain true  when $\frac{\mathbf{P}(W_n\geq x \sigma_n)}{1-\Phi \left(  x\right)}$ is replaced by $\frac{\mathbf{P}(W_n \leq -x \sigma_n)}{ \Phi \left(-  x\right)}$.

\subsection{Proof of Corollary \ref{th1.5} }\label{sec4.5}
We only need to consider the case where $\max\{\gamma_m,  \varepsilon_m,  \delta_m, m/n \}  \leq  1/10.$ Otherwise,
Corollary \ref{th1.5} holds obviously for $C$ large enough.
Denote $$\kappa_n= \alpha_0\min\{\gamma_m^{-1/4} ,\, \varepsilon_m^{-1/4},  \delta_m^{-1/4}, (m/n)^{-1/4}   \} ,$$
where $\alpha_0$ is the absolute constant given by Theorem \ref{th1}.
It is easy to see that
\begin{eqnarray}
\sup_{ x   }  \Big|\mathbf{P}( W_n  \leq x \sigma_n )  -  \Phi \left( x\right) \Big| &\leq&  \sup_{  |x|\leq \kappa_n } \Big|\mathbf{P}( W_n  \leq x \sigma_n )  -  \Phi \left( x\right) \Big| \nonumber  \\
&&  + \sup_{  |x|> \kappa_n} \Big|\mathbf{P}( W_n  \leq x \sigma_n )  -  \Phi \left( x\right) \Big|  \nonumber\\
&= & \sup_{  |x|\leq \kappa_n } \Big|\mathbf{P}( W_n  \leq x \sigma_n )  -  \Phi \left( x\right) \Big|  \nonumber\\
&&  + \sup_{   x < - \kappa_n}  \mathbf{P}( W_n  \leq x \sigma_n )  +  \sup_{   x < - \kappa_n}  \Phi \left( x\right)    \nonumber\\
&&  + \sup_{   x >   \kappa_n}  \mathbf{P}( W_n  >x \sigma_n )  +  \sup_{   x > \kappa_n} ( 1- \Phi \left( x\right) ).   \label{ineq010}
\end{eqnarray}
By Theorem \ref{th1} and the inequality $|e^x-1|\leq |x|e^{|x|},$ we have
\begin{eqnarray}
 &&\sup_{ |x|\leq \kappa_n } \Big|\mathbf{P}( W_n  \leq x \sigma_n )  -  \Phi \left( x\right) \Big| \nonumber \\
  &&\leq  \sup_{ |x|\leq \kappa_n } \Big(1-\Phi(|x|) \Big) \bigg| e^{   C_{\alpha_0}  \big( x^3 \varepsilon_m   + x^2 (  \delta_m^2   +\frac{m}n + \gamma_m|\ln \gamma_m|)  + (1+ x  ) ( \varepsilon_m   \left| \ln  \varepsilon_m
 \right|  +  \gamma_m|\ln \gamma_m|  +   \delta_m        + \sqrt{ m/n}     ) \big)}   -1 \bigg| \nonumber\\
 &&\leq  C_{\alpha_0, 1 }  \Big( \gamma_m|\ln \gamma_m|   + \varepsilon_m \left| \ln  \varepsilon_m
 \right|  + \delta_m     + \sqrt{ m/n}   \Big ). \label{ineq020}
 \end{eqnarray}
 Using the last inequality, we deduce that
 \begin{eqnarray}
\sup_{   x < - \kappa_n}  \mathbf{P}( W_n  \leq x \sigma_n )&= &  \mathbf{P}(W_n  \leq - \kappa_n  \sigma_n) \nonumber  \\
&\leq& C_{\alpha_0, 2 } \Big( \gamma_m|\ln \gamma_m|  + \varepsilon_m \left| \ln  \varepsilon_m
 \right|  + \delta_m     + \sqrt{ m/n}   \Big ) +  \Phi \left( - \kappa_n\right) \nonumber \\
&\leq& C_{\alpha_0, 3 } \Big(  \gamma_m|\ln \gamma_m|  + \varepsilon_m \left| \ln  \varepsilon_m
 \right|  + \delta_m     + \sqrt{ m/n}   \Big ).   \label{ineq030}
\end{eqnarray}
Similarly, it holds that
\begin{eqnarray}   \label{ineq040}
\sup_{   x >   \kappa_n}  \mathbf{P}(W_n   > x \sigma_n )    &  \leq&  C_{\alpha_0, 4 } \Big( \gamma_m|\ln \gamma_m|  + \varepsilon_m \left| \ln  \varepsilon_m
 \right|  + \delta_m     + \sqrt{ m/n}   \Big ).
\end{eqnarray}
It is obvious that
\begin{eqnarray} \label{ineq050}
  \sup_{   x > \kappa_n} ( 1- \Phi \left( x\right) ) =\sup_{   x < - \kappa_n}  \Phi \left( x\right)   =  \Phi \left( - \kappa_n\right)   \leq  C_{\alpha_0, 5} \Big( \gamma_m|\ln \gamma_m|  + \varepsilon_m \left| \ln  \varepsilon_m
 \right|  + \delta_m     + \sqrt{ m/n}   \Big ).
\end{eqnarray}
Combining the inequalities (\ref{ineq010})-(\ref{ineq050}) together, we obtain the desired inequality.

\subsection{Proof of Corollary \ref{co0} }\label{sec4.3}
Let $m=\sqrt{ a_n \sqrt{n}   }$. Then it holds that $m \rightarrow \infty$ as $n \rightarrow \infty.$ Thus $\max\{\gamma_m, \delta_m   \} \rightarrow 0$ as $n \rightarrow \infty .$

First, we prove that
\begin{eqnarray}\label{dfgkmsf}
 \limsup_{n\rightarrow \infty} a_n^2 \ln \mathbf{P}\bigg(  a_n W_n\in B  \bigg) \leq  - \inf_{x \in \overline{B}}\frac{x^2}{2 \sigma^2 }.
\end{eqnarray}
For any given Borel set $B\subset \mathbf{R},$ let $x_0=\inf_{x\in B} |x|.$ Then, it is obvious that $x_0\geq\inf_{x\in \overline{B}} |x|.$
Therefore, by Theorem \ref{th1},
\begin{eqnarray*}
 && \mathbf{P}\bigg(  a_n W_n \in B \bigg) \\
 &&\leq  \mathbf{P}\bigg(\, \Big|\frac{W_n}{\sigma} \Big|  \geq  \frac{ x_0}{a_n \sigma_n}\bigg)\\
 &&\leq  2\bigg( 1-\Phi \Big( \frac{ x_0}{a_n \sigma_n}\Big)\bigg) \exp\bigg\{C \bigg( (\frac{ x_0}{a_n \sigma_n})^3 \varepsilon_m + (\frac{ x_0}{a_n \sigma_n})^2 \big(\delta_m^2 +  \frac m n  + \gamma_m|\ln \gamma_m| \big) \\
  && \  \  \ \ \ \    \  \  \ \ \ \    \  \  \ \ \ \    \  \  \ \ \ \    \  \  \  \  \  \ \ \ \    \  \  \ \   +\, (1+ \frac{ x_0}{a_n \sigma_n})\big( \varepsilon_m \left| \ln  \varepsilon_m
 \right|  +\gamma_m|\ln \gamma_m|+ \delta_m  + \sqrt{\frac m n}\big)\bigg) \bigg\}.
\end{eqnarray*}
Notice that
$$   \varepsilon_m /a_n = \|X_0\|_\infty /  \sqrt{m}    \rightarrow 0$$ as $n\rightarrow \infty.$
Using   (\ref{fgsgj1}) and the fact $\lim _{n\rightarrow \infty} \sigma^2_n = \sigma^2 $,
we deduce that
\begin{eqnarray*}
\limsup_{n\rightarrow \infty} a_n^2 \ln \mathbf{P}\bigg( a_n W_n  \in B  \bigg)
 \ \leq \  -\frac{x_0^2}{2\sigma^2} \ \leq \  - \inf_{x \in \overline{B}}\frac{x^2}{2\sigma^2} ,
\end{eqnarray*}
which gives (\ref{dfgkmsf}).

Next, we prove that
\begin{eqnarray}\label{dfgk02}
\liminf_{n\rightarrow \infty} a_n^2 \ln \mathbf{P}\bigg(  a_n W_n \in B  \bigg) \geq   - \inf_{x \in B^o}\frac{x^2}{ 2 \sigma^2 } .
\end{eqnarray}
We may assume that $B^o \neq \emptyset,$ otherwise the last inequality  holds obviously because the infimum of a function over an empty set is interpreted as $\infty.$
For any $\varepsilon_1>0,$ there exists an $x_0 \in B^o,$ such that
\begin{eqnarray}
 0< \frac{x_0^2}{2 \sigma^2} \leq   \inf_{x \in B^o}\frac{x^2}{2 \sigma^2} +\varepsilon_1.
\end{eqnarray}
Without loss of generality, we may assume that $x_0>0.$
For $x_0 \in B^o,$ there exists small $\varepsilon_2 \in (0, x_0),$ such that $(x_0-\varepsilon_2, x_0+\varepsilon_2]  \subset B.$
Then it is obvious that $x_0\geq\inf_{x\in \overline{B}} x.$
By Theorem \ref{th1}, we deduce that
\begin{eqnarray*}
\mathbf{P}\bigg(  a_n W_n \in B  \bigg)   &\geq&   \mathbf{P}\bigg(  W_n  \in ( a_n^{-1}  ( x_0-\varepsilon_2), a_n^{-1} ( x_0+\varepsilon_2)] \bigg)\\
&\geq&   \mathbf{P}\Big(  W_n   > a_n^{-1}( x_0-\varepsilon_2)   \Big)-\mathbf{P}\Big( W_n  >   a_n^{-1}( x_0+\varepsilon_2) \Big)
\end{eqnarray*}
Using Theorem \ref{th1}, (\ref{fgsgj1}) and the fact $\lim _{n\rightarrow \infty} \sigma^2_n = \sigma^2 $ again, it follows that
\begin{eqnarray*}
\liminf_{n\rightarrow \infty} a_n^2\ln \mathbf{P}\bigg(a_n W_n \in B  \bigg)  \geq  -  \frac{1}{2\sigma^2}( x_0-\varepsilon_2)^2 . \label{ffhms}
\end{eqnarray*}
Letting $\varepsilon_2\rightarrow 0,$  we get
\begin{eqnarray*}
\liminf_{n\rightarrow \infty}a_n^2\ln \mathbf{P}\bigg(a_n W_n   \in B \bigg) &\geq& - \frac{x_0^2}{2\sigma^2}  \  \geq \   -\inf_{x \in B^o}\frac{x^2}{2\sigma^2} -\varepsilon_1.
\end{eqnarray*}
Because $\varepsilon_1$ can be arbitrarily small, we obtain (\ref{dfgk02}). This completes the proof of Corollary \ref{co0}.

\subsection{Proof of Theorem \ref{lemma9} }
Recall the notations in the proof of Theorem \ref{th1}. It is easy to see that
\[
\| D_{i,m}/ (n^{1/2} \sigma_n) \|_\infty \leq 2 \varepsilon_m
\]
and
\begin{eqnarray*}
\Big|\!\Big| \frac1{n\sigma_n^2} \sum_{i=1}^{k+1 }\mathbf{E}[ D_{i,m}^2 | \mathcal{F}_{(i-1)m}  ] -1 \Big|\!\Big|_\infty
 &\leq& \Big|\!\Big| \frac1{n\sigma_n^2} \sum_{i=1}^{k  }\mathbf{E}[ D_{i,m}^2 | \mathcal{F}_{(i-1)m}  ] - 1 \Big|\!\Big|_\infty +  \Big|\!\Big| \frac1{n\sigma_n^2} \mathbf{E}[ D_{k+1,m}^2 | \mathcal{F}_{km}  ]  \Big|\!\Big|_\infty \\
 &\leq &  \delta_m^2 +\frac{m}n   + 4   \varepsilon_m^2 =  \tau_m^2.
\end{eqnarray*}
Applying Lemma \ref{lemma4.2} to $\xi_i= D_{i,m}/(\sigma_n n^{1/2}  ),$  we have for all $x\geq 0,$
\begin{eqnarray}
  \mathbf{P}\bigg( W_n - \frac1{\sqrt{n}} \sum_{i=1}^{k+1 } \mathbf{E}[X_{i,m} | \mathcal{F}_{(i-1)m} ]  \geq x \sigma_n  \bigg)
   \leq   \exp\left\{-\frac{x^2}{2(1+ \tau_m^2   +\frac{2 }{3  } x\varepsilon_m )}\right\}. \nonumber
\end{eqnarray}
By an argument similar to the proof of (\ref{ineq6}), we obtain  for all $x\geq0,$
\begin{eqnarray}
 \mathbf{P} \bigg(  \Big|  \sum_{i=1}^{k+1 } \mathbf{E}[X_{i,m} | \mathcal{F}_{(i-1)m} ] \Big| \geq  x \sigma_n n^{1/2} \bigg)
 \leq
 4 \sqrt{e} \exp  \bigg\{ - \frac{x^2   }{ 2\cdot(81)^2    \gamma_m^2   } \bigg\}.
\end{eqnarray}
Using (\ref{ineq7}) again, we obtain the desired inequality.

\subsection{Proof of Proposition \ref{th3.0} }
 For each integer $n\geq1$,  let
\[
F_n(x)=\mathbf{P}(\widehat{W}_n \leq x), \ \ x \in \mathbf{R},
\]
be the cumulative distribution function of $\widehat{W}_n.$ Then  its  quantile function is define
by
\[
H_n(s)=\inf\{ x : F_n(x)\geq s\}, \ \ s \in (0, 1).
\]
Let $Z$ be a standard normal random variable. Denote
\begin{eqnarray}\label{fsffs}
Y_n = H_n( \Phi(Z)) .
\end{eqnarray}
Then  $ Y_n =_d \widehat{W}_n;$
 see Mason and Zhou \cite{MZ12}.
Denote
\begin{eqnarray}
K_n=  n^{1/2} \Big(\gamma_m|\ln \gamma_m|+ \varepsilon_m \left| \ln  \varepsilon_m
 \right|  + \delta_m +    \sqrt{\frac{m}n}  \, \Big).
 \end{eqnarray}
 By Theorem \ref{th1}, there exist an absolute constants $\beta \in (0, 1]$  and $C_\beta\geq1$ such that when $n$ is large enough, we have
for all  $ 0 \leq x \leq \beta \,  n^{1/2} \sigma_n /(m\|X_0\|_\infty ),$
\begin{eqnarray}\label{gonineq01}
\ln \Bigg| \frac { \mathbf{P}\big( Y_n>x \big) } { 1-\Phi \left( x\right)} \Bigg|  \leq     C_\beta (1+x^3) \frac{ K_n }{ n^{1/2}}
\end{eqnarray}
and
\begin{eqnarray}\label{gonineq02}
\ln \Bigg|   \frac { \mathbf{P}\big( Y_n<-x \big) } { \Phi \left(- x\right)}\Bigg|  \leq   C_\beta (1+x^3) \frac{ K_n}{ n^{1/2} }  ,
\end{eqnarray}
where  $C_{ \beta}$ depends only on  $\beta.$
By Theorem 1 of Mason and Zhou \cite{MZ12},  then whenever $n \geq   64 C_\beta^2 K_n^2$ and
\begin{eqnarray}
 |Y_n | &\leq& \Big(  \frac{\beta  \sigma_n}{m \|X_0\|_\infty}  \wedge  \frac{1}{8C_\beta K_n }  \Big)n^{1/2 }  \\
 &  \leq& \Big( \beta   \wedge  \frac{1}{8C_\beta }  \Big) \Big(\gamma_m|\ln \gamma_m|+ \varepsilon_m \left| \ln  \varepsilon_m
 \right|  + \delta_m +    \sqrt{\frac{m}n}  \,\Big)^{-1},
\end{eqnarray}
 we have
\begin{eqnarray}
  |Y_n -Z|\leq   2 C_\beta \Big( Y_n^2\, +1  \Big) \Big (\gamma_m|\ln \gamma_m|+ \varepsilon_m \left| \ln  \varepsilon_m
 \right|  + \delta_m +    \sqrt{\frac{m}n}  \, \Big),
\end{eqnarray}
which gives (\ref{fg235}) with $\alpha=\beta    \wedge  \frac{1}{8C_\beta }  $ and $C_\alpha= C_\beta$.  Notice that there exists an integer $n_0$ such that  $n \geq   64 C_\beta^2 K_n^2$ for all $n\geq n_0.$

Next we give the proof of (\ref{fg2sfs5}). Set for brevity
\[
\varsigma_n=  \gamma_m|\ln \gamma_m|+ \varepsilon_m \left| \ln  \varepsilon_m
 \right|  + \delta_m +\sqrt{\frac{m}n}.
 \]
By (\ref{fg235}), we have for all $0 \leq x \leq \frac1{4C_\alpha} \varsigma_n^{-2},  $
\begin{eqnarray}
  \mathbf{P}\Big(    |Y_n-Z| > x \, \varsigma_n \Big)
  &\leq& \mathbf{P}\Big(   |Y_n-Z|   > x \, \varsigma_n ,   |Y_n |\leq \alpha \, \varsigma_n^{-1}  \Big)
 + \, \mathbf{P}\Big(  |Y_n |> \alpha \, \varsigma_n^{-1} \Big) \nonumber \\
 &\leq& \mathbf{P}\Big( 2 C_\alpha \big( Y_n^2\, +1  \big) > x  \Big) + \, \mathbf{P}\Big(  |Y_n |> \alpha \, \varsigma_n^{-1} \Big) , \label{fsf01}
\end{eqnarray}
Notice that
\[
1-\Phi \left( x\right) \leq   \exp\{-x^2/2\}, \ \ \ \  \ x \geq 0.
\]
When $ 0 \leq x \leq \frac1{8C_\alpha} \varsigma_n^{-2} , $ by the inequalities (\ref{gonineq01}) and (\ref{gonineq02}), it holds that
\begin{eqnarray}\label{fsf02}
\mathbf{P}\Big(    2 C_\alpha \big( Y_n^2\, +1  \big) > x \Big)   &\leq&  2 \exp \bigg\{ -\frac14(  \frac{x}{2C_\alpha}     -1)  \bigg\} \nonumber\\
&\leq&    \exp \bigg\{1 -\frac x{8C_\alpha } \bigg\},
\end{eqnarray}
and that
\begin{eqnarray}\label{fsf03}
\mathbf{P}\Big(   |Y_n | > \alpha  \varsigma_n^{-1} \Big)  &\leq&   2\exp \bigg\{ -\frac 1{4} ( \alpha \varsigma_n^{-1})^2  \bigg\}\nonumber \\
&\leq&  2 \exp \bigg\{  - 2 C_\alpha \alpha^2  x \bigg\}.
\end{eqnarray}
Returning to (\ref{fsf01}), we obtain  for all $ 0 \leq x \leq \frac1{8C_\alpha} \varsigma_n^{-2} ,$
\begin{eqnarray} \label{fsf1032}
\mathbf{P}\Big(    |Y_n-Z| > x \varsigma_n \Big) &\leq&  2 \exp \Big\{ 1- c' x \Big\} ,
\end{eqnarray}
where $c'= \min\{ \frac1{8C_\alpha }  , 2C_\alpha   \alpha^2   \}.$
For $x >0,$ it is easy to see that
\begin{eqnarray}\label{fssfh023}
 \mathbf{P}\Big(   |Y_n-Z| > x \varsigma_n \Big)& \leq &  \mathbf{P}\Big(  |Y_n|  > \frac1 2 x\varsigma_n  \Big)  +  \mathbf{P}\Big(    |Z| > \frac 12 x \varsigma_n  \Big).
\end{eqnarray}
Clearly, it holds  for all $x >\frac1{8C_\alpha} \varsigma_n^{-2},$
\begin{eqnarray}
\mathbf{P}\Big(   |Z| > \frac 12 x\varsigma_n  \Big)    \leq   2 \exp \bigg\{- \frac18 x^2 \varsigma_n^2   \bigg\}
 \leq     2 \exp  \bigg\{ - \frac{1}{64 C_\alpha} x  \bigg\}. \nonumber
\end{eqnarray}
By Theorem \ref{lemma9},  there exists a positive constant $ \lambda $   such that for all $x >\frac1{8C_\alpha} \varsigma_n^{-2},$
\begin{eqnarray}
 \mathbf{P}\Big(    |Y_n|  > \frac 12 x \varsigma_n \Big)    \leq  ( 1+ 4 \sqrt{e}  )  \exp  \Big\{  - \lambda x  \Big\}  . \nonumber
\end{eqnarray}
Returning to (\ref{fssfh023}), we have for all $x >\frac1{8C_\alpha} \varsigma_n^{-2},$
\begin{eqnarray} \label{fghm53}
 \mathbf{P}\bigg(   |Y_n-Z| > x \varsigma_n   \bigg) \leq  ( 3+ 4 \sqrt{e}  ) \exp  \Big\{  - c'' x  \Big\},
\end{eqnarray}
where $c'' = \min\{ \lambda , \frac{1}{64 C_\alpha} \}.$ Combining (\ref{fsf1032}) and (\ref{fghm53}), we get the desired inequality.


\selectlanguage{english}

\end{document}